\documentclass[preprint,12pt]{elsarticle}

\usepackage{epsfig}
\usepackage{bm}
\usepackage{multirow}
\usepackage{amsmath}
\usepackage{amssymb}
\usepackage{amsfonts,graphicx,epsf}

\usepackage{lipsum}
\usepackage{epstopdf}
\usepackage{algorithmic,algorithm}
\usepackage{hyperref}

\begin{document}

\begin{frontmatter}	
	\title{Global optimization-based dimer method for finding saddle points}

	\author[yu]{Bing Yu}
	\ead{yubing93@pku.edu.cn}
	
	\author[zhang]{Lei Zhang\corref{cor1}}
	\ead{zhangl@math.pku.edu.cn}
	
	\cortext[cor1]{Corresponding author}
	
	\address[yu]{School of Mathematical Sciences, Peking University, Beijing, China.}
	\address[zhang]{Beijing International Center for Mathematical Research, Center for Quantitative Biology, Peking University, Beijing, China.}

\begin{abstract}
	Searching saddle points on the potential energy surface is a challenging problem in the rare event. When there exist multiple saddle points, sampling different initial guesses are needed in most dimer-type methods in order to find distinct saddle points. In this paper, we present a novel global optimization-based dimer method (GOD) to efficiently search saddle points by coupling ant colony optimization (ACO) algorithm with optimization-based shrinking dimer (OSD) method. In particular, we apply OSD method as a local search algorithm for saddle points and construct a pheromone function in ACO to update the global population. By applying a two-dimensional example and a benchmark problem of seven-atom island on the (111) surface of an FCC crystal, we demonstrate that GOD shows a significant improvement in computational efficiency compared with OSD method. Our algorithm offers a new framework to open up possibilities of adopting other global optimization methods to search saddle points.
\end{abstract}


\end{frontmatter}

\section{\label{sec:intro}Introduction}
   Saddle point search on the potential energy surface has attracted great attention over last two decades. The index-1 saddle point, which is often known as the transition state, is a critical point where the Hessian has one and only one negative eigenvalue. Because of the nature of its unstability, computing transition state on the potential energy surface has been a challenging problem. Meanwhile, many practical problems require the knowledge of transition states, such as calculating transition rates in chemical reactions \cite{baker1986algorithm}, predicting critical nucleus and transition pathway in phase transformation \cite{zhang2007morphology, zhang2010diffuse, zrsd2016, han2019transition}, etc.
         
   Abundant studies have focused on development of numerical methods of saddle point search. If both initial and final states are available, path-finding methods are often used to compute the minimum energy path. The most notable examples are the string method \cite{weinan2002string, erv07, dz09, ren2013climbing} and the nudged elastic band method \cite{henkelman2000improved}. When only initial state is available without knowledge of the final state, surface walking methods utilize single point associated with local information such as gradient or Hessian to search saddle points. The representative methods include the dimer-type methods \cite{henkelman1999a, ks08,pb08,zhang2012shrinking, zhang2012constrained, Xiao2014Basin, zhang2016optimization-based}, the gentlest ascent method \cite{weinan2011gentlest}, the activation-relaxation technique \cite{mousseau1998traveling,machado2011optimized}, and so on \cite{pozun2012optimizing, duncan2014biased}. 
   
   To improve the computational efficiency, optimization algorithms are often applied to speed up the saddle point search \cite{heyden2005efficient, olsen2004comparison, sheppard2008optimization}. For instance, Gao {\it et al.} proposed the iterative minimization formulation to iteratively solve a sequence of minimization problems by constructing a new objective function at each step \cite{gao2016iterative} and Gu {\it et al.} introduced the convex splitting
   method to minimize the auxiliary functional \cite{gu2018convex}. Zhang {\it et al.} proposed the optimization-based shrinking dimer (OSD) method by constructing the rotation and translation steps in the classical dimer method under an optimization framework, and then applied the Barzilai-Borwein gradient method to implement the OSD method to achieve superlinear convergence \cite{zhang2016optimization-based}. However, when there exist multiple saddle points, sampling different initial guesses are needed in order to find distinct saddle points. Thus, it remains unclear whether we can use global optimization algorithms to search unstable saddle points instead of stable minima on the potential energy surface. 
   
   Most global optimization methods can be divided into two classes: deterministic methods and stochastic methods. Deterministic global optimization finds the global solutions of an optimization problem while providing theoretical guarantees that the solution is indeed the global one within some predefined tolerance \cite{floudas2013deterministic}. Stochastic global optimization methods are iterative algorithms that generate a new candidate set of solutions from a given population using a stochastic operations, such as monte-carlo sampling \cite{li1987monte} and stochastic tunneling \cite{wenzel1999stochastic}. In particular, heuristic or metaheuristic approaches are often designed to explore the search space in order to find (near) optimal solutions with incomplete information or limited computation capacity, for instance, simulated annealing \cite{goffe1994global}, evolutionary algorithms \cite{andrzej2006evolutionary}, and swarm-based optimization algorithms \cite{poli2007particle}.

   Ant colony optimization (ACO) is a population-based metaheuristic that applies a probabilistic technique for solving optimization problems. In the classical ACO algorithm, a set of points called artificial ants are used for finding the optimal path on a weighted graph \cite{dorigo1997ant}. Artificial ants incrementally build solutions by moving on the graph. The solution construction process is stochastic and is biased by a pheromone function, which is constructed through a set of parameters associated with graph components whose values are updated by the ants at runtime.  ACO has been adapted to solve not only combinatorial problems, but also continuous optimization problems by shifting from using a discrete probability distribution to using a continuous probability density function \cite{socha2008ant}. Many efforts have been made on developing ACO variants that modify the pheromone update or combine ACO with other algorithms to improve the computational performance \cite{secckiner2013ant,liao2014ant}. 
      
   In this paper, we develop a new global optimization-based dimer (GOD) method by combining modified ACO algorithm with OSD method to search the index-1 saddle points. We first construct a pheromone function based on both the isopotential curvature and the energy gradient so that there is a high-pheromone region around saddle point. Next, we sample a group of points on the potential energy surface as the initial population in ACO. Then we use OSD as a local search method to update individual points, and apply population update in ACO algorithm in order to select points with high pheromone in a local range to continue search. During the iterations, most redundant points can be eliminated by ACO algorithm and the survivals converge to saddle points by OSD method. 
   
   The rest of this paper is organized as follows. We first briefly review OSD method in Section \ref{sec:osd} and ACO algorithm in Section \ref{sec:aco}. Then we present the construction of the pheromone function and a full description of GOD algorithm in Section \ref{sec:god}. Two numerical examples, including a two-dimensional problem and a benchmark problem of seven-atom island on the (111) surface of an FCC crystal, are shown in Section \ref{sec:examples} to demonstrate the efficiency of the proposed method. Final conclusion is given in Section \ref{sec:conclusions}.

\section{\label{sec:osd}OSD method}
   Given a potential energy function $E$ on a Hilbert space $\mathcal{H}$ of the system, we define $\nabla E(x)$ as the gradient of $E$ and choose a pair of points $x_1$ and $x_2$ as a dimer with the dimer length $l=\|x_1-x_2\|$. The dimer orientation is then given by a unit vector $v$ so that $x_1-x_2=lv$ and the (rotating) center of the dimer is defined as its geometric center, i.e., $x_{c}=(x_1+x_2)/2$. For notation convenience, let
   \begin{equation}\label{forces}
      F^i=-\nabla E(x_i),\quad i=1,2,
   \end{equation}
   be natural forces at the two endpoints of the dimer, and let
   \begin{equation}\label{force_midpoint}
      F=(F^1+F^2)/2
   \end{equation}
   be the approximated natural force at the dimer rotation center.
   
   The classical dimer method \cite{henkelman1999a} follows a two-step procedure: the dimer rotation step (the dimer is rotated toward the minimum energy configuration to obtain the dimer orientation $v$)  and the dimer translation step (the dimer climbs up by the energy ascent direction characterized by $v$). The OSD method \cite{zhang2016optimization-based} formulates dimer rotation and dimer translation to the corresponding optimization problems as follows.
  
\subsection{Dimer rotation}
  For the dimer rotation, we need to compute the ascent direction $\nu$ that is the eigenvector corresponding to the smallest eigenvalue of the Hessian, which can be translated to solve the following optimization problem:
  \begin{equation}\label{minmode}
     \min_{\nu} \frac{\nu^T H \nu}{\nu^T\nu},
  \end{equation}
  where $H$ denotes the Hessian matrix of $E$.

  The gradient function for Eq.~(\ref{minmode}) is defined as
  \begin{equation}\label{gradient function}
     g(\nu)=\frac{2}{\nu^T\nu}\left(H\nu-\lambda\nu\right),
  \end{equation}
  where $\lambda$ is the Rayleigh quotient of $(H,\nu)$, i.e. $\displaystyle\lambda=\frac{\nu^TH\nu}{\nu^T\nu}$. Denote $r(\nu)=H\nu-\lambda\nu$, which is parallel to $g(\nu)$. To avoid the direct computation of the Hessian, $\displaystyle\frac{F^2-F^1}{l}$ is used in the dimer system to approximate the action of the Hessian at the dimer center along the direction $\nu$.

  The simplest way for solving Eq.~(\ref{minmode}) is to use the steepest descent (SD) method:
  \begin{eqnarray}\label{rotation_sd}
     \nu_{k+1} & = & \nu_k-\gamma_k r(\nu_k)\nonumber\\
     & = & \nu_k-\gamma_k\left[\frac{F^2_k-F^1_k}{l}-\nu_k\nu_k^T\frac{F^2_k-F^1_k}{l}\right],
  \end{eqnarray} 
  where $\gamma_k$ is a step size. In practice, $\nu_k$ needs to be normalized at each step to guarantee $\nu_k^T\nu_k=1$. 

  Another approach is to apply the conjugate gradient (CG) method on Eq.~(\ref{minmode}). The search direction of CG is in the two-dimensional subspace composed by the gradient of the current point and the search direction of the previous step, i.e. $span\{r(\nu_k),p_{k-1}(:=\nu_k-\nu_{k-1})\}$. So the solution $\nu_{k+1}$ has the form:
  \begin{eqnarray}\label{nu_update}
     \nu_{k+1} & = & \nu_k+\eta_k^r r(\nu_k) +\eta_k^p p_{k-1}\nonumber \\
     & = &\eta_k^1 \nu_{k-1} +\eta_k^2 \nu_k+\eta_k^3 r(\nu_k).  
  \end{eqnarray} 
  Set
  \begin{eqnarray}\label{hatAB}
     \hat{A} & = & [\nu_{k-1},\nu_k,r(\nu_k)]^T H [\nu_{k-1},\nu_k,r(\nu_k)],\\
     \hat{B} & = & [\nu_{k-1},\nu_k,r(\nu_k)]^T [\nu_{k-1},\nu_k,r(\nu_k)],
  \end{eqnarray}
  where $\hat{A},\hat{B}\in \mathbb{R}^{3\times 3}$. The corresponding problem is equivalent to a three-dimensional generalized eigenvalue problem: 
  \begin{equation}\label{geigenvalue}
     \frac{\nu_{k+1}^T H \nu_{k+1}}{\nu_{k+1}^T\nu_{k+1}}=\frac{\eta^T\hat{A}\eta}{\eta^T\hat{B}\eta},
  \end{equation}
  where $\eta=[\eta^1,\eta^2,\eta^3]^T$, which is easy to calculate. If $r(\nu_k),\nu_k,\nu_{k-1}$ are collinear, we redefine $\hat{A},\hat{B}$ to be $\mathbb{R}^{2\times 2}$ matrix to make sure that $\hat{B}$ is positive. Compared to SD, CG needs to save more vectors and requires the same amount of force evaluations as SD. $H\nu_k$ is the dominated computational cost in each step for both methods. 

\subsection{Dimer translation}
  
  Since the index-1 saddle point is a maximum along the lowest curvature mode $\nu$ and a minimum along all other modes, dimer translation can be achieved by solving the following minimax optimization problem:
  \begin{equation}\label{minimax}
    \min_{x_{V^\perp}\in V^\perp}\max_{x_V\in V} E(x_V,x_{V^\perp}),
  \end{equation}
  where $V= span\{\nu\}, V^\perp=\mathcal{H}/V$.

  The gradient method for solving (\ref{minimax}) is given by
  \begin{equation}\label{dimer_translation}
    x_{k+1}=x_k+\beta_k (I-2\nu_k\nu_k^T)F_k,
  \end{equation}
  where $\beta_k$ is a step size and $F_k$ is the approximated natural force defined in Eq.~(\ref{force_midpoint}) at the $k$-th step. Thus, $-\nu_k \nu_k^T F_k$ represents an ascent direction in $V$ and $(I-\nu_k \nu_k^T) F_k$ is a descent direction in $V^\perp$.
 
  In the numerical implementation of OSD, the Barzilai-Borwein (BB) gradient method has been successfully applied to achieve a superlinear convergence. The step sizes $\gamma_k$ in Eq.~(\ref{rotation_sd}) and $\beta_k$ in Eq.~(\ref{dimer_translation}) are chosen by BB step sizes \cite{zhang2016optimization-based}. We shrink the dimer length $l(t)$ to approach a small constant after the dimer rotation and translation steps in order to guarantee the local convergence \cite{zhang2012shrinking}.

\section{\label{sec:aco}Modified ACO algorithm}
  ACO is a metaheuristic method for solving global optimization problems with a randomly given population (points set). Unlike the classical ACO algorithm on graphs, we proposed a modified ACO algorithm by using two steps: local search and population update.
  
\subsection{Step 1: local search}
  The step of local search is to drive the individual points in the current population to approach the local solutions by applying some dynamics or optimization algorithms. We denote the local search function by $LocalSearch(\cdot)$:
  \begin{equation}\label{localsearch}
     x^{*} \leftarrow LocalSearch(x).
  \end{equation}
  For instance, the steepest descent method, conjugate gradient method and Quasi-Newton methods can be often used as the local search methods in the minimization problems. Once a point converges to the local minimum by following local search, we move it from the population set to the solution set.
  
\subsection{Step 2: population update}
  The step of population update is to update the population by keeping the high-pheromone points and deleting the low-pheromone points in probability. The pheromone function $ph(\cdot)$ is an evaluation function that increases the pheromone value when close to the extrema and reaches the local or global maximum at the extrema. For instance, a simple pheromone function can be chosen as $ph(x)=-E(x)$ when searching the minima of the energy function.

 The $PopulationUpdate(\cdot)$ algorithm is as follows: 
  
  \begin{algorithm}[H]
  	\caption{PopulationUpdate}
  	\label{alg:p_update}
  	\begin{algorithmic}[1]
  		\REQUIRE ~~\\
  		   Population $P:=\{x_1,x_2,\cdots,x_n\}$, solution set $S$, parameters $\delta_1\gg\delta_2>0$. 
  		\ENSURE ~~\\
  		   Updated population $\hat{P}$.
  		\STATE{$\hat{P}\leftarrow\emptyset$ $\;\%$initialize output} 
  		\STATE{$P\leftarrow \{x_i\in P|distance(x_i,S)>\delta_2\} $ $\;\%$delete redundant points close to the solution set}
  		\FOR {$x_i\in P$}
  		\STATE{$P_i\leftarrow\{x_j|distance(x_i,x_j)<\delta_1,x_j\in P\}$ $\;\%$define a neighborhood of $x_i$}
  		\IF{$ph(x_i)=max\{ph(x_j)|x_j\in P_i\}$}
  		\STATE{$\hat{P}\leftarrow \hat{P}\bigcup \{x_i\}$ $\;\%$select the point with the highest pheromone}
  		\ELSE
  		\STATE{$x_j=RouletteWheelSelection(P_i,ph(\cdot))$}
  		\STATE{$\hat{P}\leftarrow \hat{P}\bigcup \{x_j\}$ $\;\%$select the other points in probability}
  		\ENDIF
  		\ENDFOR
  		\RETURN $\hat{P}$.
  	\end{algorithmic}
  \end{algorithm}  
  
  In Algorithm \ref{alg:p_update}, we delete the redundant points close to the current solution set $S$ to avoid finding the repeated saddle points because the non-degenerate stationary points are isolated according to the Morse lemma \cite{milnor2016morse}.
    
  $RouletteWheelSelection(\cdot)$ is a genetic operator to randomly select points, which is also known as the fitness proportionate selection \cite{lipowski2012roulette-wheel}. For the subpopulation $P_i$, the probability of $x_j\in P_i$ being selected is 
  \begin{equation}\label{RWS}
     prob(x_j)=\frac{ph(x_j)}{\sum_{x_{j'}\in P_i}ph(x_{j'})}
  \end{equation}  
  in $RouletteWheelSelection(\cdot)$. Therefore, the points with high pheromone have larger probability to be selected than the points with low pheromone. 

\section{\label{sec:god}Framework of GOD}
  
\subsection{Pheromone function}
 
To search saddle points, it is critical to choose a suitable pheromone function for population update in ACO algorithm. Notice that both saddle points and local extrema are the stationary points that satisfy $\nabla E(x)=0$, we need to impose the other geometric property besides the gradient in the pheromone function in order to distinguish saddle points and extrema. The eigenvector corresponding to the smallest eigenvalue of Hessian matrix is used for the ascent direction in the gentlest ascent method \cite{weinan2011gentlest}. Moreover, the isopotential curvature $\kappa$ is a description of how fast the contour tangent direction changes as moving along the isopotential contour. It has been used to control the dimer dynamics for searching saddle points \cite{Xiao2014Basin} and is defined as 
  \begin{equation}\label{isopotential_curvature}
    \kappa(x) = - \min_{c\perp \nabla E(x)}\frac{c^T H(x) c}{\|\nabla E(x)\|}.
  \end{equation}
Here $c$ is a unit vector that satisfies $c\perp \nabla E(x)$. Denoting $\lambda_1\le\lambda_2\le\cdots\le\lambda_n$ as the eigenvalues of $H$ and $\nu_1$ as $\min_{c\perp \nabla E}(c^T H c)$, then we have $\lambda_1\le\nu_1\le\lambda_2$. In the neighbourhood of an index-1 saddle point, the eigenvalues satisfy $\lambda_1<0<\lambda_2$, indicating there is a curve on which $\kappa=0$ extends from the index-1 saddle point. 

Since the isopotential curvature is unbounded around critical points, we redefine $\kappa$ as
   \begin{equation}\label{kappa}
      \kappa(x) = - \min_{c\perp \nabla E(x)}\frac{c^T H(x) c}{c^T c},
   \end{equation}  
 which is a constrained Rayleigh quotient. 
   
   By combining the gradient and the isopotential curvature, the pheromone function $ph(\nabla E(x),\kappa(x))$ is now defined as:
   \begin{equation}\label{pheromone}
      ph(\nabla E(x),\kappa(x))=\frac{\alpha}{1+a|\kappa(x)|}+\frac{1-\alpha}{1+b\|\nabla E(x)\|},
   \end{equation}
   where $\alpha\in (0,1)$ is a parameter to represent the $\kappa$'s importance to the gradient and $a,b$ are some positive parameters to adjust the norm of $\kappa$ and gradient.   
   
   As in Eq.~(\ref{kappa}), the computation of $\kappa$ can be translated to a constrained optimization problem such as Eq.~(\ref{minmode}):
   \begin{equation}\label{kappa2}
     \min_{(c,g)=0} \frac{c^THc}{c^Tc},
   \end{equation}
   where $g$ denotes the gradient of energy at $x$. The CG method in Eq.~(\ref{nu_update}) is used to solve Eq.~(\ref{kappa2}), and $\nu_k$ in OSD can be used as the initial guess. The only difference is that we need to project initial point and search direction into the subspace $span\{g\}^\perp$.

\begin{figure}[htbp]
	 \centering
	 \includegraphics[width=1\textwidth]{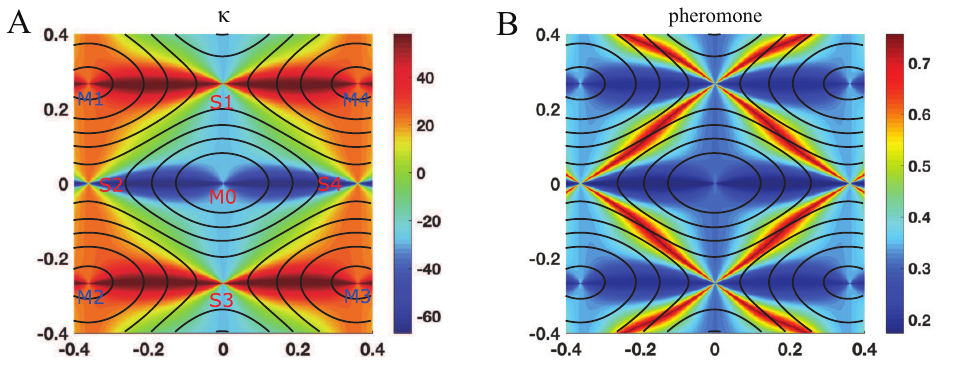}
	 \caption{An example of $B_2$ function that has one minimum M0, four maxima M1-M4, and four index-1 saddle points S1-S4. (A): modified $\kappa$ in Eq.~\eqref{kappa}. (B): pheromone function in Eq.~(\ref{pheromone}). The black lines are the contours of the $B_2$ function. The parameters used in Eq.~(\ref{pheromone}) are $\alpha=0.75,a=0.05,b=100$.}
	 \label{fig:k_ph}
   \end{figure}

 As an illustration of  $\kappa$ in Eq.~(\ref{kappa}) and the pheromone function in Eq.~(\ref{pheromone}), we apply a two-dimensional (2D) $B_2$ function \cite{socha2008ant}:
   \begin{eqnarray}\label{2d_func}
      f(x,y) & = & x^2+2y^2-0.3\cos(3\pi x)\nonumber\\
      & & -0.4\cos(4\pi y)+0.7, x,y\in[-0.4,0.4].
   \end{eqnarray}
   In the Fig.~\ref{fig:k_ph}(A), $\kappa$ is negative around the minimum and positive around the maxima, and the boundary $\kappa=0$ intersects the separatrix at all saddle points. By taking advantage of the information on both $\kappa=0$ and $\|\nabla E(x)\|=0$, there are high-pheromone regions near saddle points (Fig.~\ref{fig:k_ph}(B)), which can provide the guide of suitable paths for searching saddle points.

\subsection{GOD Algorithm}
   By substituting $LocalSearch(\cdot)$ by OSD method and applying the pheromone function $ph(\cdot)$ in population update, we can achieve a practical GOD algorithm as below:
   \begin{algorithm}[H]
   	\caption{GOD algorithm}
   	\label{alg:GOD}
   	\begin{algorithmic}[1]
   		\REQUIRE ~~\\
   		   Initial population $P^1:=\{x_1^1,x_2^1,\cdots,x_n^1\}$, parameters $\delta_1\gg\delta_2>0$, number of local search in each iteration $step_{ls}$, and maximum of iteration $step_{max}$.
   		\ENSURE ~~\\
   		   Saddle points set $S$.
   		\STATE{Initialize: $S:=\emptyset$.}
   		\FOR {$k=1$ to $step_{max}$}
   		  \STATE{$\tilde{P} \leftarrow\emptyset$ $\quad\%$ initialize temporary set}
   		  \FOR{$x_i^k\in P^k$} 	 	     
   		     \STATE{$x_i^{k*} \leftarrow OSD(x_i^k,step_{ls})$  $\quad\%$ do $step_{ls}$-step local search}
   		     \IF{$x_i^{k*}$ converges}
   		        \STATE{$S\leftarrow S\cup\{x_i^{k*}\}$ $\quad\%$ put saddle point into solution set} 
   		     \ELSE
   		        \STATE{$\tilde{P}\leftarrow \tilde{P}\cup\{x_i^{k*}\}$ $\quad\%$ put other points into temporary set}
   		     \ENDIF
   		  \ENDFOR	
   		  \IF{$\tilde{P}=\emptyset$}
   		     \STATE{\textbf{exit}$\;\%$ stop when population is empty}
   		  \ENDIF
   		  \STATE{$P^{k+1}\leftarrow PopulationUpdate(\tilde{P},S,\delta_1,\delta_2)$ $\;\%$ update the population}   		  
   		\ENDFOR
   		\RETURN $S$
   	\end{algorithmic}
   \end{algorithm}

In practice, there is unnecessary to update the pheromone and the population in every iteration because calculating $\kappa$ needs force evaluations as many as that in dimer rotation. Thus, we impose a parameter $step_{ls}$ that represents the number of local search steps in each iteration. Moreover, if the population does not change (i.e., $P^{k+1}=\tilde{P}$) after population update, indicating most remaining points are well separated, we can stop the population update and only apply local search to continue saddle point search till all points converge. Since the ACO has no influence on the searching directions of $LocalSearch(\cdot)$, the convergence of the algorithm can be guaranteed by OSD method \cite{zhang2016optimization-based}. 
  
 \section{\label{sec:examples}Numerical examples}

\subsection{Two-dimensional example}
\begin{figure}[htbp]
	\centering
	\includegraphics[width=1\textwidth]{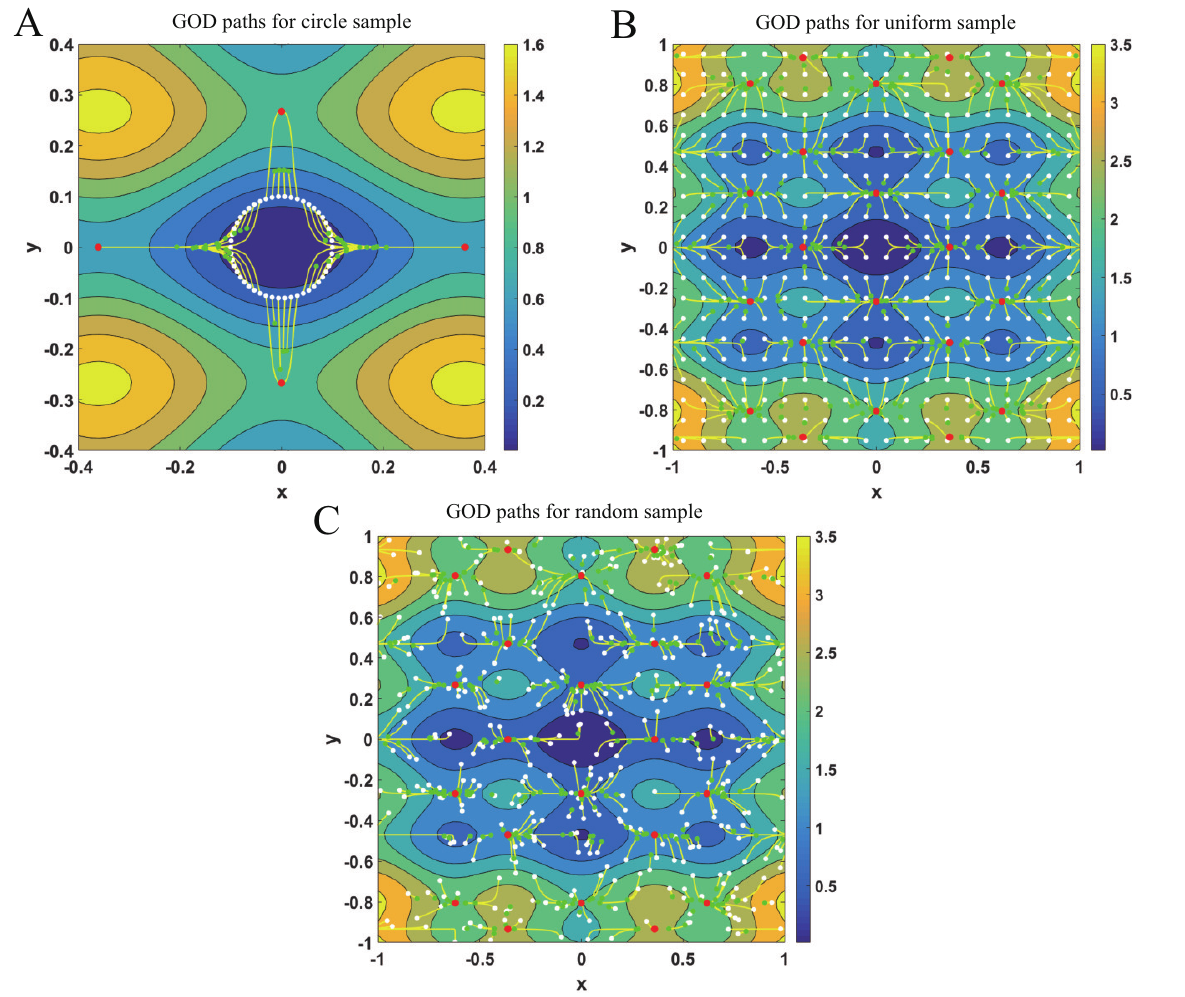}
	\caption{2D example with three different initial samplings: (A) circle sampling; (B) uniform sampling; (C) random sampling. White dots represent the initial points, green dots are the points deleted in the halfway, and red dots are the saddle points. The yellow lines represent the dynamic paths by GOD.}
	\label{fig:GOD_paths}
\end{figure}

  We first apply the 2D $B_2$ function to test the efficiency of GOD algorithm, which is given by 
   \begin{equation*}\label{2d_func_2}
      f(x,y) =  x^2+2y^2-0.3\cos(3\pi x) -0.4\cos(4\pi y)+0.7.
   \end{equation*}  
 In Fig.~\ref{fig:GOD_paths}(A), we sample 50 initial points in a circle (white dots) around the local minimum $(0,0)$ in the domain $[-0.4,0.4]\times [-0.4, 0.4]$. By using GOD algorithm, most points are deleted in the halfway (labeled by green dots) and four saddles (red dots) can be found when the GOD algorithm converges.

We also test this example in a larger domain $[-1,1]\times [-1, 1]$, in which there are totally $22$ saddle points. We use 400 initial points by applying two different sampling: one is uniform sampling and the other is random sampling. As shown in Fig.~\ref{fig:GOD_paths}(B) and (C), all saddles can be found and the most points are deleted along dynamical paths to avoid redundant computation for the same saddles.

We count the number of force evaluations of GOD algorithm and compare it with OSD method in Table~\ref{2d_sample}. We set $step_{ls}=1$ for simplicity, i.e. update population at every $LocalSearch(\cdot)$ step in GOD. If only OSD method is applied, each initial point will converge to a saddle, resulting a large number of force evaluations. GOD algorithm can reduce huge computation cost and require less force evaluations ($<20\%$) than OSD mthod from all three different samplings while finding all saddles successfully.

\begin{table}[htbp]
	\centering
	\begin{tabular}{|c|c|c|c|c|c|}
		\hline
		sampling & method & initial points & force evaluations & saddles \\ 
		\hline
		circle & OSD & 50 &  7270 & 4 \\
		circle & GOD & 50 &  1249 & 4 \\
		\hline
		uniform & OSD & 400  & 44212 & 22 \\
		uniform & GOD & 400 &  8716 &  22 \\
		\hline
		random & OSD & 400  & 45581 & 22 \\
		random & GOD  & 400 & 8352 & 22 \\
		\hline  
	\end{tabular}
		\caption{\label{2d_sample} Comparison of the force evaluations for the 2D example by applying OSD and GOD methods.}	
\end{table}

\subsection{Seven-atom island on the (111) surface of an FCC crystal}
  Next, we implement a high-dimensional example by putting the seven-atom island on the (111) surface of an FCC crystal in three dimensions (3D) \cite{zhang2016optimization-based}, which often serves as a benchmark problem for searching saddle points. The cluster consists of seven atoms with the substrate made by a 6-layer slab, each of which contains $56$ atoms. The bottom three layers in the slab are fixed and the other $56\times 3+7=175$ atoms are free to move. The atoms interact via a pairwise additive Morse potential
  \begin{equation}\label{morse_potential}
     V(r)=D(e^{-2\alpha(r-r_0)}-2e^{-\alpha(r-r_0)}).
  \end{equation}
Parameter setting is $D=0.7102$, $\alpha=1.6047$, and $r_0=2.8970$. 

\begin{figure}[H]
	\centering
    \includegraphics[width=0.8\textwidth]{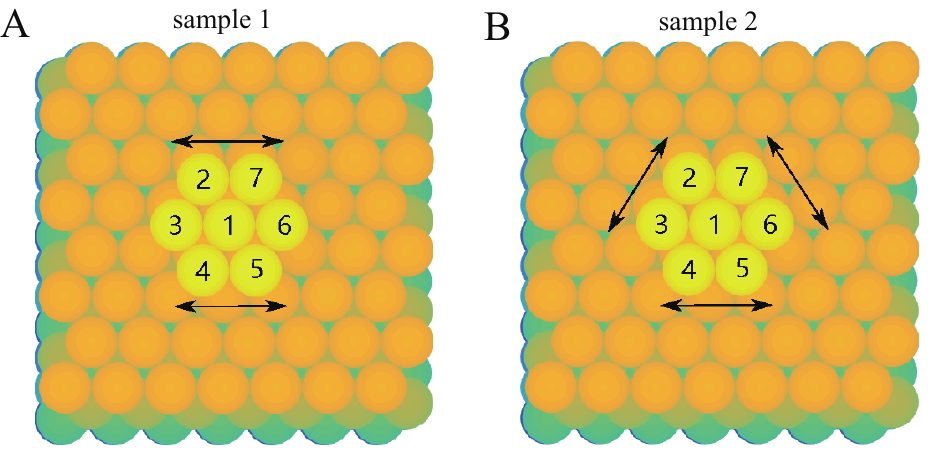}   
	\caption{Two initial samplings for the seven-atom island: (A) sample 1: shift the parallel $(2,7)$ and $(4,5)$ atoms simultaneously to get $100$ initial points; (B) sample 2: shift parallel atoms along six directions to get $300$ initial points.}
	\label{shift}
\end{figure}

  We first sample two different initial populations and use the OSD method to search possible saddle points. The cputime and the number of force evaluations will be taken as the reference in order to test the GOD algorithm. Since it is clear that the density of samplings can influence the efficiency of GOD, we choose two initial samplings by shifting the top seven atoms in Fig.~\ref{shift} to avoid too dense or too sparse sampling. Sample 1 is a group of $100$ points, which shifts the $(2,7)$ and $(4,5)$ atoms along the directions of arrows in Fig.~\ref{shift}(A). In sample 2, we shift the parallel atoms along six directions in Fig.~\ref{shift}(B) to get $300$ initial points.

\begin{figure}[htbp]
	\centering
	\includegraphics[width=1\textwidth]{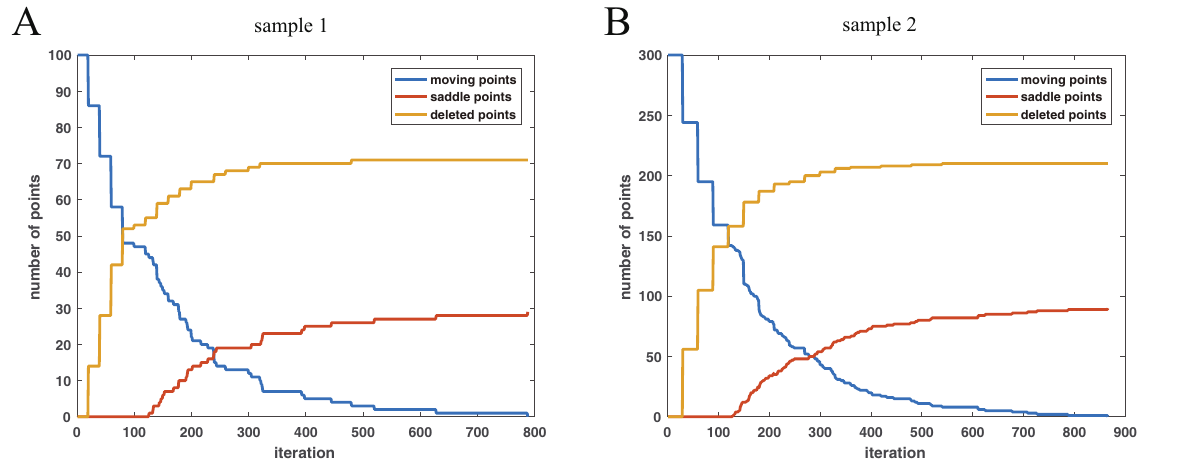}
	\caption{Seven-atom island example with two different initial samplings: (A) sample 1; (B) sample 2. Blue line represents the number of moving points, yellow line represents the number of deleted points, and red line corresponds to the number of saddle points. Some parameters are $\delta_1=0.7, a=b=10, \alpha=0.5$. $step_{ls}=20$ for sample 1 and $step_{ls}=30$ for sample 2.}
	\label{fig:antsvsiteration}
\end{figure}

Fig.~\ref{fig:antsvsiteration} shows the numerical results for both sample 1 and sample 2. 

Blue line represents the number of moving points for searching saddle points, yellow line shows the number of deleted points, and the red line corresponds to the number of saddles we have found. During the iterations, the number of the moving points quickly decreases and some jumps appear because of the population update step. The stopped points turn into either the deleted points or the saddle points till the last moving point converges to the saddle point. 

\begin{table}[htbp]
	\centering
	\begin{tabular}{|c|c|c|c|c|}
		\hline
		sampling & $step_{ls}$ & cputime(s) & force evaluations & saddles \\ 
		\hline
		1 & OSD & 325.781  & 80946   & 29 \\
		\hline
		1 & 5 & 134.313  & 31101 &  20 \\
		1 & 10                                  & 160.875  & 37713 &  25 \\
		1 & 15                                  & 173.641  & 41324 &  28 \\
		1 & 20                                  & 188.375  & 45072 &  29 \\
		\hline  
		2 & OSD  & 1022.720  & 250941   & 92 \\
		\hline
		2 & 10   & 493.563   & 116206   & 78 \\
		2 & 20   & 543.813   & 130202   & 87 \\
		2 & 30   & 618.906   & 140624   & 92 \\
		\hline  
	\end{tabular}
		\caption{\label{label_step_ls}Effect of $step_{ls}$ in GOD algorithm. Some parameters are $\delta_1=1.0$, $a=b=10,\alpha=0.5$.}	
\end{table}

There exists a trade-off between local search and population update. Population update aims to delete close points to save the computational cost but may lose some points that will converge to different saddles by local search because of different searching directions. Thus, we can tune the parameter $step_{ls}$ to balance local search steps and population update. Table \ref{label_step_ls} shows the smaller $step_{ls}$ results in the less cputime and force evaluations, but with a cost of losing some saddle points compared with OSD method. Increasing $step_{ls}$ can eventually recover all saddle points while reducing $>40\%$ computational cost, indicating the optimal $step_{ls}$ can enhance the efficiency of GOD algorithm.

Besides the parameter $step_{ls}$, we apply the sample 1 to test the sensitivity of  the neighborhood parameter $\delta_1$ in population update and the parameters in the pheromone function in Table \ref{label_parameter}. Default values of parameters here are $step_{ls}=15,\delta_1=0.75, \alpha=0.5, a=b=10$. Since $\delta_1$ determines the influence region of pheromone function in population update, smaller $\delta_1$ costs more computational cost and larger $\delta_1$ may result in the lose of saddle points. On the other hand, the change of the parameters $\alpha, a, b$ in the pheromone function does little influence on the computational cost and saddles, indicating GOD algorithm is robust to these parameters.

\begin{table}[htbp]
	\centering
	\begin{tabular}{|c|c|c|c|c|}
		\hline
		\multicolumn{2}{|c|}{parameter} & cputime(s)  & force evaluation & saddles \\ 
		\hline
		\multicolumn{2}{|c|}{OSD} &  325.781  & 80946   & 29  \\
		\hline
		\multirow{4}{*}{$\delta_1$} & 
		  1.0  & 173.641  & 41324 &  28 \\
		& 0.75 & 202.875 & 48479  &  29 \\
		& 0.5  & 228.984  & 53341 &  29 \\
		\hline
		\multirow{3}{*}{$\alpha$} & 
		  0.25 & 203.359 & 47994  & 29 \\
		& 0.5  & 202.875 & 48479  & 29 \\
		& 0.75 & 194.859 & 46359  & 29 \\ 
		\hline
		\multirow{4}{*}{$a$} & 
		  1    & 207.125 & 49098  & 28 \\
		& 10   & 202.875 & 48479  & 29 \\
		& 100  & 201.266 & 47826  & 29 \\
		& 1000 & 195.484 & 46662  & 28 \\
		\hline
		\multirow{4}{*}{$b$} & 
		  1    & 207.547 & 49434  & 29 \\
		& 10   & 202.875 & 48479  & 29 \\
		& 100  & 205.406 & 48826  & 28 \\
		& 1000 & 205.109 & 48892  & 29 \\
		\hline  
	\end{tabular}
	\caption{\label{label_parameter}Effect of $\delta_1$ in population update and the parameters in the pheromone function for Sample 1.}
\end{table}

\section{\label{sec:conclusions}Conclusions}
  
  In this paper, we propose a novel GOD algorithm for searching index-1 saddle points on the potential energy surface. This algorithm is a combination of modified ACO algorithm with OSD method. More specifically, we apply OSD method as a local search method for saddle points and construct a pheromone function in ACO algorithm for global population update. As demonstrated by a 2D example and the seven-atom island problem, GOD algorithm can greatly reduce the cputime and force evaluations, and shows a significant improvement in computational efficiency compared with OSD method.
  
  One secret ingredient of GOD algorithm is the construction of the pheromone function that guides the population distribution on the potential energy surface, which could be further improved to enhance the performance of GOD algorithm. Moreover, other surface walking methods can be naturally adopted in the framework of GOD. In the cases of high-dimensional energy functions, how to do efficient sampling is still a challenging problem and remains unsolved. Some attempts have been made for the minimum mode following method \cite{plasencia2016improved} to reduce the number of force evaluations and avoid the repeated saddle points. Overall, good flexibility of GOD algorithm provides a framework to develop efficient global optimization algorithms for searching saddle points and solve application problems.
    
\section*{Acknowledgments}
  We would like to thank Professors Qiang Du and Weiqing Ren for fruitful discussions. This work was supported by the National Natural Science Foundation of China
  11622102, 11861130351, 11421110001.

  \bibliography{references}

\end{document}